\documentclass[a4,article,prl,nofootinbib, twocolumn]{revtex4-1}

\usepackage{graphicx}
\usepackage{dcolumn}
\usepackage{bm}
\usepackage{color}
\usepackage{amsmath}%
\usepackage{amsfonts}%
\usepackage{amssymb}%
\usepackage{breqn}
\usepackage{fixfoot}

\makeatletter
\let\cat@comma@active\@empty
\makeatother

\begin{document}
\title{Up-Hill Diffusion Creating Density Gradient - What is the Proper Entropy?}
\author{N. Sato and Z. Yoshida}
\affiliation{Graduate School of Frontier Sciences, The University of Tokyo,
Kashiwa, Chiba 277-8561, Japan}
\date{\today}

\begin{abstract}
It is always some constraint that yields any nontrivial structure from statistical averages.
As epitomized by the Boltzmann distribution, the energy conservation is often the principal constraint
acting on mechanical systems. Here, we investigate a different type: the topological constraint
imposed on `space'. 
Such constraint emerges from the null space of the Poisson operator linking energy gradient to phase space velocity, and appears as an adiabatic invariant altering the preserved phase space volume at the core
of statistical mechanics.
The correct measure of entropy, built on the distorted invariant measure, behaves consistently with the second law of thermodynamics. The opposite behavior (decreasing entropy and negative entropy production) arises in arbitrary coordinates.   
An ensamble of rotating rigid bodies is worked out.
The theory is then applied to up-hill diffusion in a magnetosphere.
\end{abstract}


\maketitle

\begin{normalsize}

There are plenty of examples that seemingly violate the principle of entropy maximization. 
So-called up-hill diffusion, creating density gradients, is often observed in multi-phase fluids and solids undergoing spinodal decomposition \cite{Hashimoto,Cahn}, in metallic alloys \cite{Dar}, nanoporous materials \cite{Lau}, and  magmas \cite{Les}.
By separating the different components of the mixture, Helmholtz free energy achieves a local minimum, characterized by non-uniform concentrations, that is stable against fluctuations \cite{Cahn}.   
With a completely different mechanism, astronomical plasmas accumulate within the magnetic fields of stars and planets through the process of inward diffusion \cite{Schulz,Boxer,Yoshida2010,SatoStoc} and generate an heterogeneous density profile. That the driving force is not the energy constraint is made apparent by the experimental observation of non-neutral plasma particles climbing up the potential hill \cite{Yoshida2010}, as well as by numerical calculations concerning their thermal equilibrium \cite{SatoTherm}.
Here, the underlying principle
is the self-organization of a quasi-stationary state, governed by long-range interactions, that feeds upon the 
\textit{topological constraints} (typically in the form of adiabatic or Casimir invariants \cite{TC}) affecting canonical phase space.  
As long as the invariants are preserved, the ordered architecture, arising from the integral manifolds foliating phase space, seems to be conflicting with the second law of thermodynamics.
Once the invariants are broken, the quasi-stationary state is destroyed and the systems progressively approach thermal death.  
Accretion of galaxies under the action of gravitation \cite{Lynden-Bell,Chavanis}, ferromagnetism mediated by the magnetic field \cite{Antoniazzi,Pakter}, spontaneous creation of planetary magnetospheres through the electromagnetic interaction \cite{Schulz,Boxer,Yoshida2010}, vortical structures in magnetofluids preserving helicities \cite{YosMah2002}, living organisms harvesting `negentropy' \cite{Schrodinger}, self-organization of data flows in information theory \cite{Ikegami}   
are some of the most paradigmatic examples of such ordered structures that grow on the topological invariants affecting the relevant `phase space'.

In the present paper, we 
study the non-equilibrium statistical mechanics of Hamiltonian systems subjected to the aforementioned topological constraints. In particular, we show that the entropy defined on the invariant measure (the preserved phase space volume) of the system behaves consistently with the second law of thermodynamics.  
Due to the non-covariant nature of differential entropy \cite{Jaynes,Jaynes2}, the time evolution of the uncertainty measured in arbitrary coordinates may `flip', and appear as an entropy decrease in the Cartesian perspective.  
It is the Jacobian of the coordinate change that yields the ordered structure, while the probability distribution is flattened in the proper variables.

The theory, which finds its roots in the phenomenological observation that particle density in planetary magnetopsheres tends to be homogenized in the magnetic coordinates \cite{Bir, Has}, shows that the proper phase space upon which statistical mechanics can be formulated differs from the \textit{a priori} variables used to represent a general physical system. 
These findings may pave the way for a new and rigorous understanding of the statistical mechanics governing constrained systems.



We start with a short review of the Hamiltonian formalism.
Hamiltonian mechanics is the result of interaction between matter (energy or Hamiltonian function $H$) and space (Poisson operator $\mathcal{J}$) according to the equation:

\begin{equation}
\boldsymbol{v}=\mathcal{J}\nabla H, \label{HM}
\end{equation}

\noindent where $\boldsymbol{v}=\dot{\boldsymbol{x}}$ is the flow velocity in $n$-dimensional phase space. 
(\ref{HM}) admits two typologies of constants of motion: those that can be ascribed to the specific form of the Hamiltonian function, i.e. to the properties of matter, and the so-called Casimir invariants that originate from the eigenvectors with $0$-eigenvalue (the null space or kernel) of the Poisson operator, i.e. from the properties of space. 
This second kind of invariants, which limits the accessible regions of phase space as a result of the constraining environment, is at the core of the theory developed in the present work.
Specifically, due to antisymmetry $\mathcal{J}^{T}=-\mathcal{J}$, whenever the operator $\mathcal{J}$ has a kernel $\boldsymbol{\xi}$ such that $\mathcal{J}\boldsymbol{\xi}=\boldsymbol{0}$, the system is subjected to \textit{topological constraints}:

\begin{equation}
\boldsymbol{\xi}\cdot \boldsymbol{v}=0.\label{TC}
\end{equation}


\noindent (\ref{TC}) can be thought as the formal definition of topological constraint.
We remark that the above result holds for any Hamiltonian, and even if $\mathcal{J}$ does not satisfy the Jacobi identity (see \cite{f1}). However, thanks to Darboux's theorem \cite{Morrison1998}, the Jacobi identity ensures that the kernel is integrable, i.e. that a Casmir invariant exists:

\begin{equation}
\boldsymbol{\xi}=\lambda\nabla C,\label{ITC}
\end{equation} 

\noindent where, for now, we assumed that the rank of $\mathcal{J}$ is $n-1$ (see \cite{rank}), and the two functions $\lambda$ and $C$ are integration factor and Casimir invariant ($\dot{C}=\lambda^{-1}\boldsymbol{\xi}\cdot\boldsymbol{v}=0$) respectively.  


It is now useful to make some considerations on the non-covariant nature of differential entropy.    
Extension of Shannon's discrete entropy to continuous probability
distributions is a delicate process \cite{Jaynes,Jaynes2}. Indeed, the quantity:

\begin{equation}
\tilde{S}=-\int_{V}{p\left(\boldsymbol{x}\right)\log{p\left(\boldsymbol{x}\right)}dV}\label{STil}
\end{equation}

\noindent is not, in general, the entropy of the continuous probability distribution $p\left(\boldsymbol{x}\right)$ on the volume element $dV=dx^{1}\wedge ... \wedge dx^{n}$. The reason is that $\tilde{S}$ is not covariant, i.e. its value changes depending on the chosen coordinate system, 
and (\ref{STil}) tacitly assumes that $dV$ is an invariant measure. Unfortunately, this is not always the case and (\ref{STil}) has to be amended with Jaynes' functional:

\begin{equation}
S^{J}=-\int_{V}{p\left(\boldsymbol{x}\right)\log{\left(\frac{p\left(\boldsymbol{x}\right)}{\mathfrak{I}\left(\boldsymbol{x}\right)}\right)}dV},\label{SJ}
\end{equation}

\noindent where the Jacobian $\mathfrak{I}\left(\boldsymbol{x}\right)$ compensates the coordinate dependence of the logarithm. In the Hamiltonian picture, one can always find a time-independent function $\mathfrak{I}\left(\boldsymbol{x}\right)$ nullifying the Lie derivative of $\mathfrak{I}dV$ with respect to the dynamical flow $\boldsymbol{v}$, i.e. 
such that $\mathfrak{L}_{\boldsymbol{v}}\mathfrak{I}\left(\boldsymbol{x}\right)dV=0$.
The obtained $\mathfrak{I}$ with (\ref{SJ}) will then give the desired covariant form of entropy. 
It is useful to recast (\ref{SJ}) as below: 

\begin{equation}
\Sigma=-\int_{V_{I}}{P\left(\boldsymbol{y}\right)\log{P\left(\boldsymbol{y}\right)}dV_{I}}.\label{S}
\end{equation}

\noindent Here, $P$ is the probability distribution of $\boldsymbol{y}$ and $dV_{I}$ is the invariant measure of the system 
$dV_{I}=dy^{1}\wedge ... \wedge dy^{n}=\mathfrak{I}dV$ 
satisfying $\mathfrak{L}_{\boldsymbol{u}}dV_{I}=0$, with $\boldsymbol{u}=\dot{\boldsymbol{y}}$. 


We are now ready to test the theory with a simple $3$D example.
In $3$D, equation (\ref{HM}) can always be cast in the form $\boldsymbol{v}=\boldsymbol{w}\times\nabla H$, 
where $\boldsymbol{w}$ is a properly chosen vector (see \cite{f2}).  
The Euler's rotation equation for the motion of a rigid body with angular momentum $\boldsymbol{x}$ and moments of inertia $I_{x}$, $I_{y}$, and $I_{z}$ can be obtained by setting $H=\left(x^{2}/I_{x}+y^{2}/I_{y}+z^{2}/I_{z}\right)/2$ and $\boldsymbol{w}=\boldsymbol{x}$.
The kernel $\boldsymbol{\xi}$ associated to this operator, i.e. the topological constraint (\ref{TC}) affecting the phase space of a rigid body, is soon identified to be $\boldsymbol{\xi}=\boldsymbol{x}$. Indeed, $\boldsymbol{\xi}\cdot\boldsymbol{v}=\boldsymbol{x} \cdot\boldsymbol{x}\times\nabla H=0$.
At the same time, one can verify that the Jacobi identity (see \cite{f3})
is satisfied $\boldsymbol{x}\cdot\nabla\times\boldsymbol{x}=0$, making the system Hamiltonian. The Jacobi identity also guarantees
that the kernel is integrable (remember (\ref{ITC})) to give the integration factor $\lambda=1$ and the Casimir
invariant $C=\boldsymbol{x}^{2}/2$, so that $\boldsymbol{w}=\nabla C$. Furthermore, the invariant measure turns out to be $dV_{I}=dx\wedge dy\wedge dz$, as follows from $\nabla\cdot\boldsymbol{v}=0$. 
Since this is the original statistical measure, one can directly apply (\ref{STil}) to define the entropy of an ensemble of such rigid bodies. However, suppose that we consider a slightly more complicated rotation pattern, such as: 

%
%

\begin{equation}
\boldsymbol{v}=\lambda(\boldsymbol{x})\nabla \frac{\boldsymbol{x}^{2}}{2}\times\nabla H,\label{RB}
\end{equation}


\noindent where, for example, $\lambda=e^{z^{2}/2}$. Since $\dot{z}\propto\lambda$, high values of $z$ will be less probable and (\ref{RB}) may represent the anisotropic rotation of a rigid body that tends to spin around the axis with angular momenta $x,y$. (\ref{RB}) still satisfies the Jacobi identity, and thus represents an Hamiltonian system with the same Casimir element $C$.
However, the invariant measure becomes:

\begin{equation}
dV_{I}=e^{-z^{2}/2}dx\wedge dy\wedge dz= dC\wedge d\chi \wedge dz,\label{IMRB}
\end{equation}  

\noindent where we introduced new coordinates $(C,\chi,z)$, with $\chi=e^{-z^{2}/2}\arctan\left(y/x\right)$.
%
Separating the constant of motion $C$, the new $2$D canonical equations are:

\begin{equation}
\boldsymbol{u}=
\begin{bmatrix} \dot{\chi} \\ \dot{z} \end{bmatrix}=\begin{bmatrix}-H_{z}\\H_{\chi}\end{bmatrix}.
\label{RBY}
\end{equation}

\noindent 
One can verify that (\ref{RBY}) is divergence free. 

In order to study the statistical mechanics of the new system, we now consider an ensemble of objects obeying (\ref{RBY}) and let them interact by adding to the Hamiltonian an interaction potential $\phi$. Its ensemble average must go zero $\left\langle\phi\right\rangle=0$, since the total energy of the system has to be preserved. In addition, and this is the key point of the paper, there are grounds for the ergodic hypothesis in the novel coordinates $(C,\chi,z)$ (and not in the original variables $(x,y,z)$) because of the invariant measure (\ref{IMRB}). In other words, it is licit to exchange ensemble averages with time averages \textit{only} on (\ref{IMRB}):


\begin{equation}
\begin{split}
0= &\left\langle d\phi\right\rangle=\left\langle\phi_{\chi}\right\rangle d\chi+\left\langle\phi_{z}\right\rangle dz=\\ &\bar{\phi}_{\chi}d\chi+\bar{\phi}_{z}dz=
\bar{\Gamma}_{\chi}\left(t\right)d\chi+\bar{\Gamma}_{z}\left(t\right)dz,
\end{split}\label{EH}
\end{equation} 

%

\noindent with $\Gamma_{\chi}$ and $\Gamma_{z}$ Gaussian white noises and where the bar indicates long-time averaging. In (\ref{EH}) first we substituted ensemble averages with time averages, and then represented the various components with random processes of zero time average. We remark that this would not have been possible in the original coordinates $(x,y,z)$, as they are not measure preserving. Neglecting the constant $C$, the equations accounting for the interaction become:

\begin{equation}
\begin{bmatrix}\dot{\chi}\\\dot{z}\end{bmatrix}=\begin{bmatrix}-H_{z}-\Gamma_{z}\\H_{\chi}+\Gamma_{\chi}\end{bmatrix}\label{EoMI}.
\end{equation}

\noindent Note that, while the Hamiltonian is no more a constant, $C$ is still a Casimir invariant:
the rigid bodies will explore the surface of phase space defined by $C$.

The next step is to build the Fokker-Planck equation associated to (\ref{EoMI}). We refer the reader to \cite{Gardiner,SatoStoc} for a detailed description of the procedure. The result is:

\begin{dmath}
\frac{\partial P}{\partial t}=H_{z}\frac{\partial P}{\partial \chi}-H_{\chi}\frac{\partial P}{\partial z}+\frac{1}{2}D_{\chi}\frac{\partial^{2}P}{\partial\chi^{2}}+\frac{1}{2}D_{z}\frac{\partial^{2}P}{\partial z^{2}}.\label{FPE}
\end{dmath}

\noindent Here $P$ is the probability distribution on $\left(\chi,z\right)$ and $D_{\chi}$, $D_{z}$ are the diffusion coefficients associated with the white noises. Finally, we seek for an explicit expression of the entropy production rate $\sigma$ of the system. Define the Fokker-Planck velocity $\boldsymbol{Z}$ to be the vector field such that (\ref{FPE}) is written as $\partial_{t}P=-\nabla\cdot\left(\boldsymbol{Z}P\right)$. Then, recalling (\ref{S}):

\begin{equation}
\frac{d\Sigma}{dt}=\int_{V_{I}}{\left\{P\nabla\cdot \boldsymbol{Z}+\nabla\cdot\left[P\log\left(P\right)\boldsymbol{Z}\right]\right\}}dV_{I}.
\end{equation}

\noindent The first term represents the ensemble average of the Fokker-Planck velocity divergence, while the second factor can be cast as a surface integral representing entropy flow out $L$. It is straightforward to deduce that:

\begin{subequations}
\begin{align}
\sigma&=\left\langle\nabla\cdot \boldsymbol{Z}\right\rangle,\label{EPRZ}\\
L&=-\int_{V_{I}}\nabla\cdot\left[P\log\left(P\right)\boldsymbol{Z}\right]dV_{I}.
\end{align}
\end{subequations}

\noindent Substituting the expression of $\boldsymbol{Z}$ in (\ref{EPRZ}), we obtain:

\begin{equation}
\sigma=-\frac{1}{2}D_{\chi}\left\langle\frac{\partial^{2}\log{P}}{\partial\chi^{2}}\right\rangle-\frac{1}{2}D_{z}\left\langle\frac{\partial^{2}\log{P}}{\partial z^{2}}\right\rangle.
\end{equation}

\noindent In figure \ref{fig1} we report the results of the numerical simulation of (\ref{FPE}). 
The Entropy $\Sigma$, defined on the invariant measure (\ref{IMRB}) of the system, behaves consistently with the second law of thermodynamics and the associated entropy production $\sigma$ is positive. On the contrary, the wrong measure of entropy $\tilde{S}=-\int f\log{f}dV=\Sigma+\langle\lambda\rangle$, defined by the distribution function $f$ on the original phase space $dV=dx\wedge dy\wedge dz$, decreases. Furthermore, diffusion flattens the distribution $P$ and since preservation of particle number requires $PdV_{I}=fdV$, $f=P/\lambda$ creates an ordered structure by approaching $f\propto \lambda^{-1}$.

\begin{figure}[h]
\hspace*{-0.65cm}\centering
\includegraphics[scale=0.28]{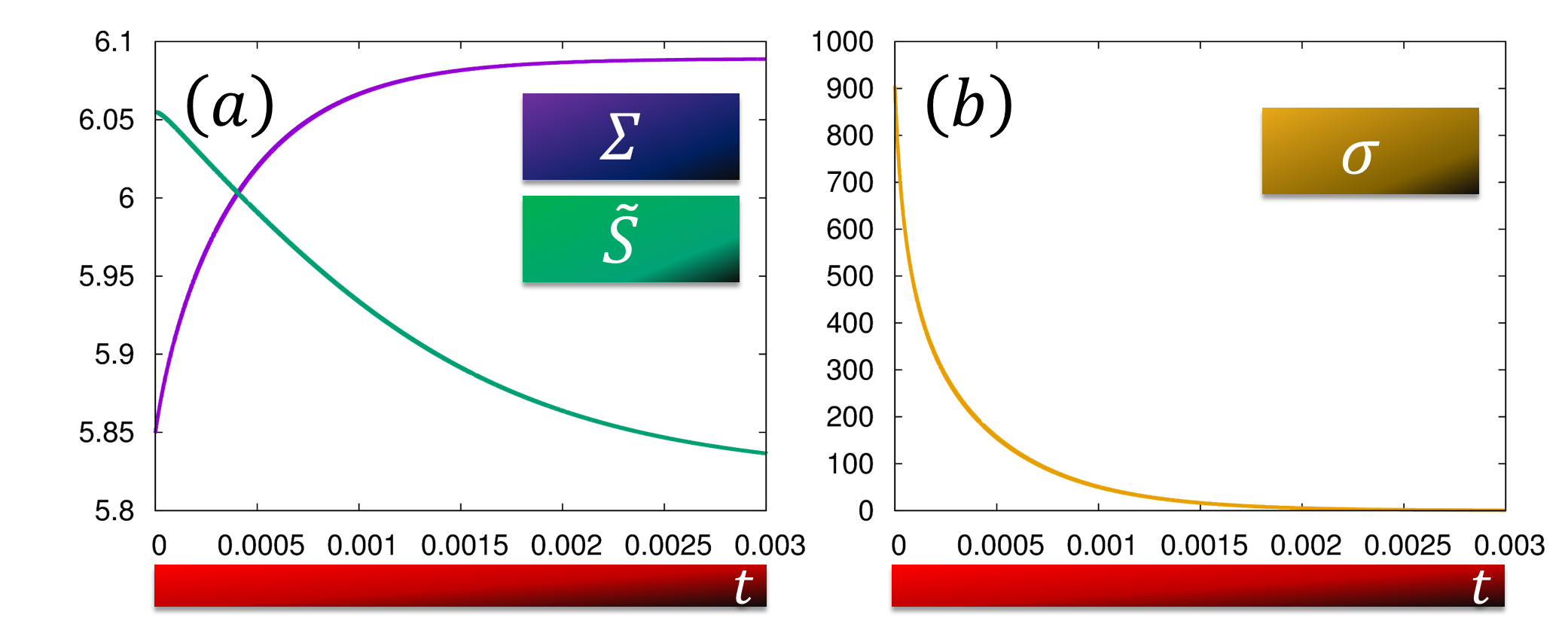}
\caption{\footnotesize (a): $\Sigma$ and $\tilde{S}$ as a function of time $t$. (b): $\sigma$ as a function of time $t$. Arbitrary units are used. Initial condition is a flat distribution $f$ on $dV$.}
\label{fig1}
\end{figure}
 

Let us show how the theory can be applied to the study of a real self-organizing system: a magnetosphere. In astronomical plasmas, charged particles are trapped by planetary magnetospheres as they spiral around the magnetic field $\boldsymbol{B}=\nabla\psi\times\nabla\theta$, where $\psi=\psi(r,z)$ is the flux-function and $\theta$ the toroidal angle of a cylindrical coordinate system $(r,z,\theta)$. This dynamics (cyclotron motion) is characterized by preservation of the magnetic moment $\mu=mv^{2}_{\perp}/2B=const$, where $m$ is the particle mass, $v_{\perp}$ the particle velocity perpendicular to magnetic field lines, and $B=\lvert\boldsymbol{B}\rvert$. 
Because of the topological constraint $\mu$, it turns out \cite{SatoStoc,SatoHeat} that the invariant measure of magnetized particles is $dV_{I}=d\mu\wedge dv_{\parallel}\wedge dl\wedge d\psi\wedge d\theta=Bd\mu\wedge dv_{\parallel}\wedge dx\wedge dy\wedge dz=BdV$, where $l$ and $v_{\parallel}$ are length and velocity along $\boldsymbol{B}$ respectively. The electromagnetic interaction diffuse the constrained particles on the statistical measure $dV_{I}$ and maximize the associated entropy $\Sigma$. Due to the inhomogeneous Jacobian $B$, the process will appear as creating density gradients and temperature anisotropy in the Cartesian perspective, while the entropy $\tilde{S}$ defined on $dV$ is minimized. This scenario is exemplified in figures \ref{fig2} and \ref{fig3} obtained by simulation of the Fokker-Planck equation derived in \cite{SatoStoc,SatoHeat}.


\begin{figure}[h]
\hspace*{-0.65cm}\centering
\includegraphics[scale=0.28]{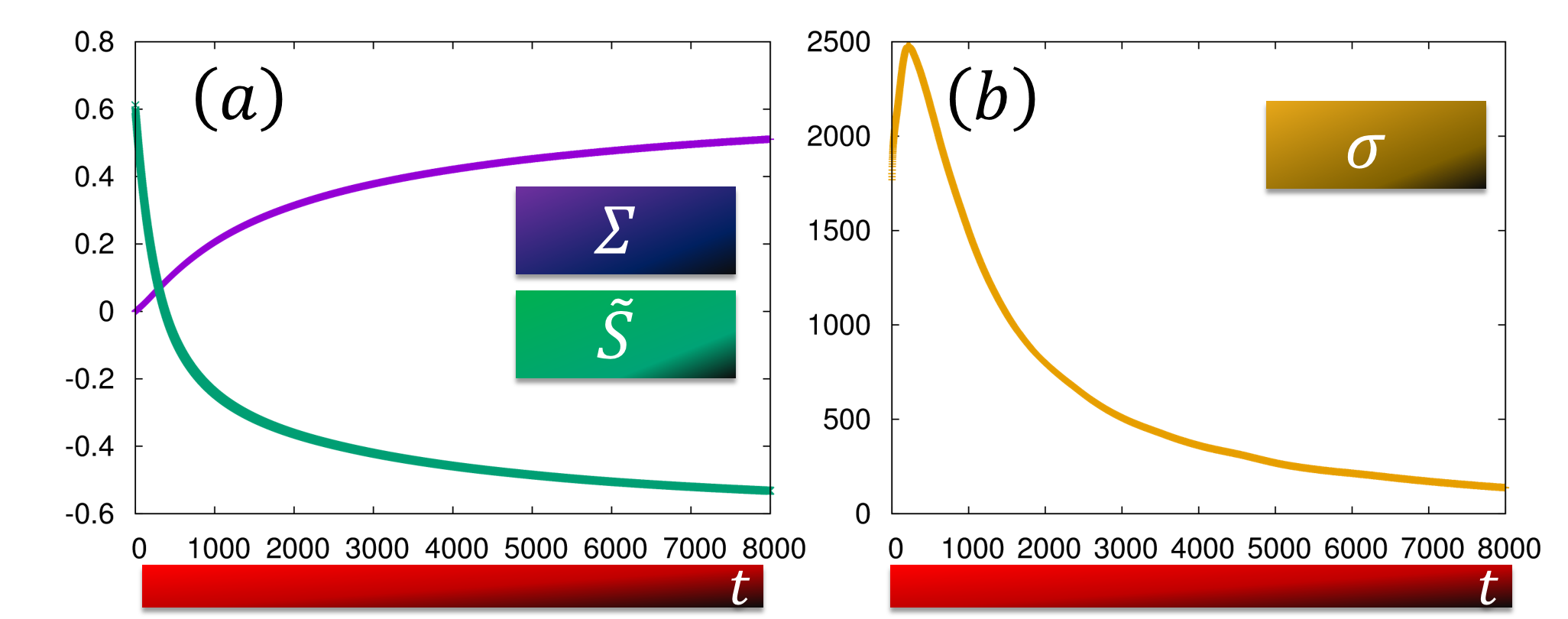}
\caption{\footnotesize (a): $\Sigma$ and $\tilde{S}$ as a function of time $t$. (b): $\sigma$ as a function of time $t$. Arbitrary units are used. Initial condition is a Maxwell-Boltzmann distribution.}
\label{fig2}
\end{figure}

\begin{figure}[h]
\hspace*{-0.7cm}\centering
\includegraphics[scale=0.5]{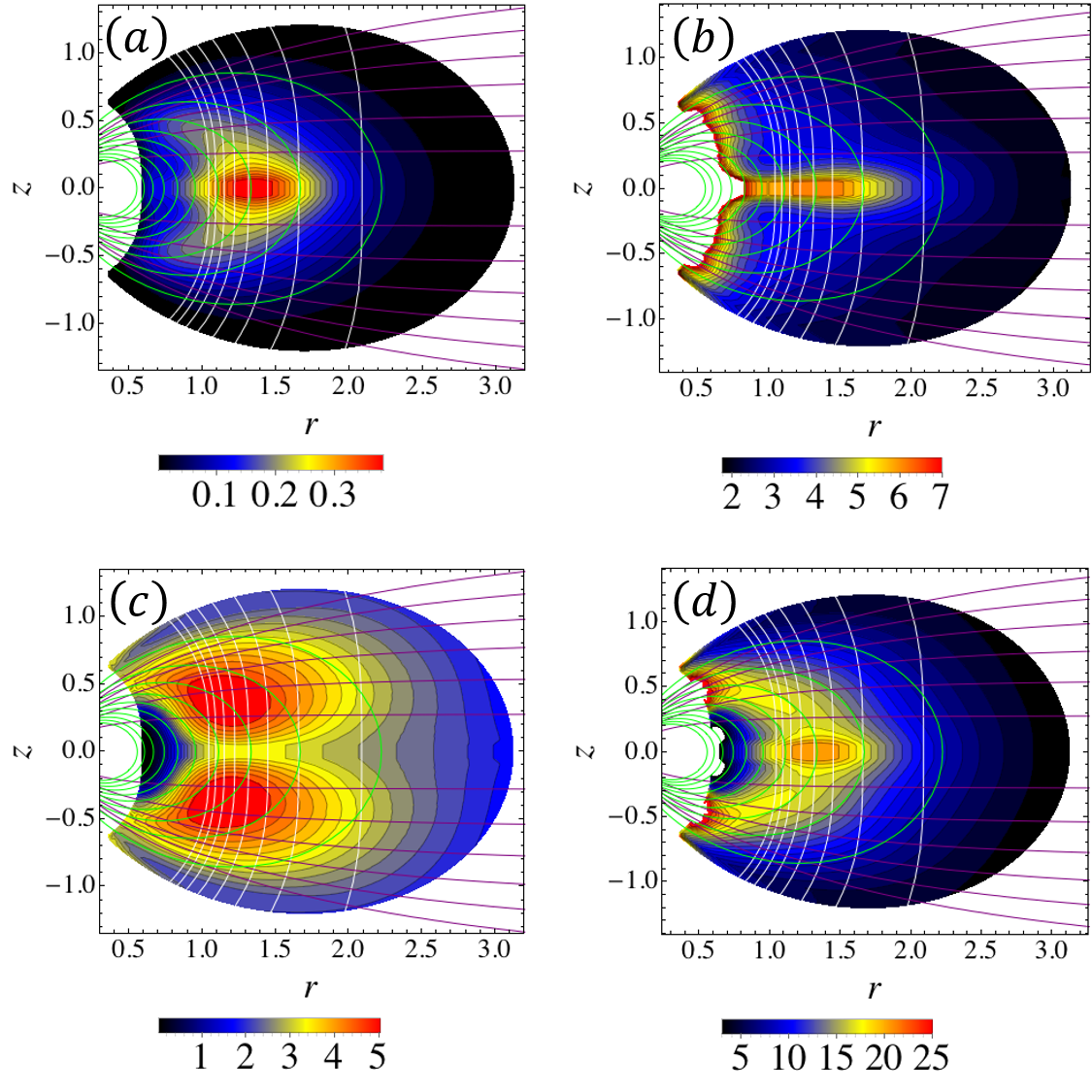} 
\caption{\footnotesize Self-organized plasma after entropy maximization. (a): spatial profile of particle density (a.u.). (b): temperature anisotropy $T_{\perp}/T_{\parallel}$. (c): parallel temperature $T_{\parallel}(eV)$. (d): perpendicular temperature $T_{\perp}(eV)$. White, green, and purple lines represent contours of $B$, $\psi$, and $l$.}
\label{fig3}
\end{figure}

This research was supported by JSPS KAKENHI Grant Nos. 23224014 and 15K13532.



\end{normalsize}


\begin{thebibliography}{99}

\bibitem{Hashimoto} T. Hashimoto, K. Matsuzaka, E. Moses, and A. Onuki, Phys. Rev. Lett. \textbf{74}, 1 (1995). 



\bibitem{Cahn} J. W. Cahn, Acta Met. \textbf{9}, 795-801 (1961).


\bibitem{Dar} L. S. Darken, Trans. AIME \textbf{180}, 430-438 (1949).


\bibitem{Lau} A. Lauerer, T. Binder, C. Chmelik, E. Miersemann, J. Haase, D. M. Ruthven, and J. Karger,  Nat. Comm. \textbf{6}, 7697 (2015).



\bibitem{Les} C. E. Lesher, J. Geophys. Res. \textbf{99}, B5 pp. 9585-9604 (1994).


\bibitem{Schulz} M. Schulz and L. J. Lanzerotti, 
\textit{Particle Diffusion in the Radiation Belts} (Springer, New York, 1974).

\bibitem{Boxer} A. C. Boxer, R. Bergmann, J. L. Ellsworth, D. T. Garnier, J. Kesner, M. E. Mauel, and P. Woskov, 
Nature Phys. \textbf{6}, 207 (2010).

\bibitem{Yoshida2010} Z. Yoshida, H. Saitoh, J. Morikawa, Y. Yano, S. Watanabe, Y. Ogawa,  
Phys. Rev. Lett \textbf{104}, 235004 (2010).

\bibitem{SatoStoc} N. Sato and Z. Yoshida, J. Phys. A: Math. Theor. \textbf{48}, 205501 (2015).


\bibitem{SatoTherm} N. Sato, N. Kasaoka, and Z. Yoshida, Phys. Plasmas \textbf{22}, (4) 042508 (2015). 


\bibitem{TC} A rigorous definition of topological constraint based on the degeneracy of the Poisson algebra will be given later on. When the constraint is integrable, it is called a Casimir invariant in the context of non-canonical Hamiltonian mechanics. The correspondence with adiabatic invariants of classical mechanics is discussed in Z. Yoshida and P. J. Morrison, in \textit{Nonlinear physical systems: spectral analysis, stability and bifurcation}, edited by O. N. Kirillov and D. E. Pelinovsky (ISTE and John Wiley and Sons, 2014), Chap. 18, pp. 401-419 and Z. Yoshida and S. M. Mahajan, Prog. Theor. Exp. Phys. \textbf{2014} 073J01 (2014).



\bibitem{Lynden-Bell} D. Lynden-Bell and R. Wood,
Mon. Not. R. Astron. Soc. \textbf{138}, 495 (1968).

\bibitem{Chavanis} P. H. Chavanis, 
\textit{Dynamics and Thermodynamics of Systems with Long-Range Interactions: An Introduction} (Springer, 2002), pp. 208-289.


\bibitem{Antoniazzi} A. Antoniazzi, D. Fanelli, S. Ruffo, and Y. Y. Yamaguchi. 
Phys. Rev. Lett. \textbf{99}, 040601 (2007).

\bibitem{Pakter} R. Pakter and Y. Levin,
Phys. Rev. Lett. \textbf{106}, 200603 (2011).



\bibitem{YosMah2002} Z. Yoshida and S. M. Mahajan, 
Phys. Rev. Lett. \textbf{88}, 095001 (2002).
 
\bibitem{Schrodinger}  E. Schrodinger, \textit{What is life- the physical aspect of the living cell}, (Cambridge University Press, 1944).


\bibitem{Ikegami} T. Ikegami and M. Oka, ICAART 2014, pp. 237-242 (2014).






\bibitem{Jaynes} E. T. Jaynes, \textit{Probability theory the logic of science} (Cambridge University Press, 2003), pp. 374-376. 


\bibitem{Jaynes2} E. T. Jaynes,  Phys. Rev.  \textbf{106}, 4 pp. 620-630 (1957). 



\bibitem{Bir} T. J. Birmingham, T. G. Northrop, and C. G. Falthammar, Phys. Fluids \textbf{10}, 11 (1967).


\bibitem{Has} A. Hasegawa, Phys. Scr. T \textbf{116}, 72-74 (2005).




\bibitem{f1} This identity is satisfied by any Poisson operator and is the essential feature of the algebraic structure of Hamiltonian systems. Defining the Poisson bracket of two functions $f$ and $g$ as $\{f,g\}=\nabla f\cdot\mathcal{J}\nabla g$, the identity for three functions $f$, $g$, and $h$ reads $\{f,\{g,h\}\}+\{g,\{h,f\}\}+\{h,\{f,g\}\}=0$.

\bibitem{Morrison1998} P. J. Morrison, Rev. Mod. Phys. \textbf{70}, 2 (1998).


\bibitem{rank} If the rank of $\mathcal{J}$ is $n-k$ and the Jacobi identity is satisfied, there will be $k$ Casimir
invariants. Again, this is a consequence of Darboux's theorem \cite{Morrison1998}.

\bibitem{f2} If $\boldsymbol{w}=\boldsymbol{B}/B^{2}$ with $\boldsymbol{B}$ the magnetic field, and setting $H=\phi$, with $\phi$ the electric potential, one obtains the equations of motion for a magnetized particle performing $\boldsymbol{E}\times\boldsymbol{B}$ drift. In vacuum, $\boldsymbol{B}=\nabla\xi$ for some pontetial $\xi$. In this case, $\xi$ is a Casimir invariant.

\bibitem{f3} In $3$D the Jacobi identity reads as $\boldsymbol{w}\cdot\nabla\times\boldsymbol{w}=0$. When satisfied, the constraint $\boldsymbol{w}\cdot \boldsymbol{v}=0$ is integrable \cite{Szekeres}, and thus two scalars $\lambda$ and $C$ can be found such that $\boldsymbol{w}=\lambda\nabla C$. Furthermore, $\mathfrak{I}=\lambda^{-1}$ defines an invariant measure since $\nabla\cdot\left(\mathfrak{I}\boldsymbol{v}\right)=0$, with the result that, in $3D$, Jacobi identity, Casimir invariant, and invariant measure imply each other provided that $\boldsymbol{v}=\boldsymbol{w}\times\nabla H$.

\bibitem{Szekeres} P. Szekeres, \textit{A Course in Modern Mathematical Physics}, (Cambridge University Press, 2004), pp. 454-455.


\bibitem {Gardiner} C. W. Gardiner, \textit{Handbook of Stochastic Methods for Physics, Chemistry and the Natural Sciences}, 2nd ed., (Springer-Verlag, 1985).



\bibitem{SatoHeat} N. Sato, Z. Yoshida, and Y. Kawazura, Plasma Fus. Res. \textbf{11}, 2401009 (2016).






\end{thebibliography}
\end{document}